\documentclass{amsart}
\usepackage{amsfonts,amssymb,amsmath,amsthm,mathrsfs}
\usepackage{url}
\usepackage[dvips]{epsfig}
\newtheorem{thrm}{Theorem}[section]

\newtheorem{prop}[thrm]{Proposition}
\newtheorem{cor}[thrm]{Corollary}

\newtheorem{example}[thrm]{Example}
\begin{document}
\author[C.~A.~Mantica and L.~G.~Molinari]
{Carlo~Alberto~Mantica and Luca~Guido~Molinari}
\address{C.~A.~Mantica: I.I.S. Lagrange, Via L. Modignani 65, 
20161, Milano,  and I.N.F.N. sezione di Milano, Via Celoria 16, 20133 Milano, Italy 
-- L.~G.~Molinari (corresponding author): 
Physics Department Aldo Pontremoli,
Universit\`a degli Studi di Milano and I.N.F.N. sezione di Milano,
Via Celoria 16, 20133 Milano, Italy.}
\email{carlo.mantica@mi.infn.it, luca.molinari@mi.infn.it}
%\subjclass[2010]{53B20, 53C50 (Primary), 83C20 (Secondary)}
%\keywords{Conformally quasi-recurrent manifold,
%Weyl compatible tensor, Petrov types, Lorentzian metric.} 
%
\title[GENERALIZED CURVATURES AND COMPATIBLE TENSORS]
{The Jordan Algebras of  Riemann, Weyl\\ and curvature compatible tensors}

\begin{abstract}
Given the Riemann, or the Weyl, or a  generalized curvature tensor $K$, a symmetric tensor $b_{ij}$ is named `compatible'
with the curvature tensor if  $b_i{}^m K_{jklm} + b_j{}^m K_{kilm} + b_k{}^m K_{ijlm}=0$.
%The symmetric tensors that are compatible with the Riemann, or the Weyl, or a 
%generalized curvature tensor $K$, 
Amongst showing known and new properties, we prove that they form a special Jordan algebra, i.e. the symmetrized 
product of K-compatible tensors is K-compatible. 
%Identities by Lovelock are recovered, and new generalized curvature tensors are obtained.
\end{abstract}
\date{26 September 2019}
\maketitle
\section{\bf Introduction} 
Let $(M, g)$ be a $n-$dimensional Riemannian or pseudo-Riemannian manifold, and $K_{jklm}$ a  generalized
curvature tensor (the Riemann, the Weyl, or any tensor with the algebraic properties of the Riemann
tensor). In ref.\cite{EXTDS} we introduced this concept: a symmetric tensor $b_{ij}$ is $K-$compatible if 
\begin{align}
b_i{}^m K_{jklm} + b_j{}^m K_{kilm} + b_k{}^m K_{ijlm}=0. \label{bK}
\end{align}
We name $(K,b)$ a {\em compatible pair}. The motivation was the following theorem \cite{EXTDS}: if $b_{ij}$ is $K-$compatible with eigenvectors $X,Y,Z$ and eigenvalues  $x,y,z$ with $z\neq x,y$, then:
\begin{align}
K_{ijlm} X^iY^jZ^m =0. \label{XYZ}
\end{align}
It extends a result by Derdzi\'nski and Shen \cite{DerdzShen} who proved the same for the Riemann
tensor, with the hypothesis that $b_{ij}$ is a Codazzi tensor, $\nabla_i b_{jk}=\nabla_j b_{ik}$. 
Despite the increased generality, the replacement of the Codazzi condition with the algebraic condition \eqref{bK}, enabled a far simpler proof of the new theorem. 

Equation \eqref{bK} with Riemann's tensor originally appeared in a paper by Roter, 
on conformally symmetric spaces (\cite{Roter1972} lemma 1).
Riemann and Weyl compatible tensors were studied in refs. \cite{RCT,WCT,Deszcz2013}. \\
Examples of Riemann compatible tensors are the Codazzi tensors \cite{EXTDS}, 
the Ricci tensors of Robertson-Walker or perfect-fluid 
generalized Ro\-bert\-son-Wal\-ker space-times \cite{PFGRW},  the second fundamental form and the Ricci tensor 
of a hypersurface embedded in a (pseudo)Riemannian manifold \cite{WCT}, the Ricci tensors of `weakly Z-symmetric' manifolds 
($\nabla_i Z_{jk} = A_i Z_{jk} + B_j Z_{ik} + D_k Z_{ij}$  with $Z_{ij}=R_{ij}+\varphi g_{ij}$, $A_k-B_k$ closed 1-form) \cite{WZSYM} that include  `weakly Ricci-symmetric' ones ($\varphi =0$) \cite{WRSYM} and others (see \cite{GENPSRS,Chaki}), or `pseudosymmetric manifolds' \cite{PSEUDOSYM} 
($[\nabla_i,\nabla_j] R_{klmp} = L Q_{klmpij}$, where $L\neq -1/3$ is a
scalar function and $Q$ is the Tachibana tensor built with the Riemann and Ricci tensors).\\
A Riemann compatible tensor is also Weyl compatible, but not the opposite. The Ricci tensors of 
G\"odel (\cite{Deszcz2014}, th.2), or pseudo-Z symmetric space times \cite{Mant_Suh3} are Weyl compatible.

In sections 2 and 3 we review Riemann and Weyl compatible tensors, with some new results and 
examples, and their relation with known identities by Lovelock. Then, in sections 4, 5 and 6, we investigate the algebraic properties of 
generalized curvature tensors and $K-$compatible tensors. The main result is that the latter form a 
{\em special Jordan algebra}, i.e. the set of $K-$compatible tensors is closed for the symmetrized product.

\section{ Riemann  compatible tensors}
A symmetric tensor is Riemann compatible if: % \cite{EXTDS,RCT}:
\begin{align}
b_i{}^m R_{jklm} +b_j{}^m R_{kilm} + b_k{}^m R_{ijlm}=0.  \label{bR}
\end{align}
The relation may be written $b_{(i}{}^m R_{jk)lm}=0$, where $(ijk)$ denotes the sum on cyclic permutations of the indices. Contraction 
with the metric tensor $g^{jl}$ gives $R_{km}b^m_i  - b_k{}^m R_{mi}=0$ i.e. $b$ commutes with the Ricci tensor. Contraction
with $b^{jl}$ gives $b_i{}^m R_{jklm}b^{jl}  + b_k{}^m R_{ijlm}b^{jl} =0 $ i.e. $b$ commutes with the symmetric tensor 
$\hat R_{jm} = R_{jklm} b^{kl}$.

\begin{example} 
Codazzi tensors are Riemann compatible.\\
Proof: in the identity $[\nabla_i,\nabla_j]b_{kl} = -R_{ijl}{}^m b_{km} -R_{ijk}{}^m b_{ml}$ sum on
cyclic permutations of $ijk$. The first Bianchi identity $R_{(ijk)}{}^m=0$, gives:
\begin{align*}
[\nabla_i, \nabla_j]b_{kl} +[\nabla_j, \nabla_k]b_{il}+[\nabla_k, \nabla_i]b_{jl}=
-(b_i{}^m R_{jklm} +b_j{}^m R_{kilm} + b_k{}^m R_{ijlm}). 
\end{align*} 
The left hand side is zero for Codazzi tensors. 
\end{example}
\begin{example}
If  $\nabla_j A_k = p_j A_k$, then $A_iA_j$ is Riemann compatible.\\ 
Proof:  $A_i [\nabla_j,\nabla_k]A_l = A_i(\nabla_j p_k - \nabla_k p_j)A_l = A_l[\nabla_j,\nabla_k]A_i$.
Then $A_i R_{jkl}{}^m A_m = A_l R_{jki}{}^m A_m$;
the sum on cyclic permutations of $ijk$ gives zero in r.h.s.  
\end{example}
\subsection{Codazzi deviation} In ref.\cite{RCT} we introduced the natural concept of 
{\em Codazzi deviation} of a symmetric tensor:
\begin{align}
\mathscr C_{jkl}=\nabla_j b_{kl}-\nabla_k b_{jl}. 
\end{align}
\noindent
Properties: $\mathscr C_{jkl}=-\mathscr C_{kjl}$, $\mathscr C_{jkl}+\mathscr C_{klj}+\mathscr C_{ljk} =0$, and
\begin{align}
\nabla_i\mathscr C_{jkl}+\nabla_j\mathscr C_{kil}+\nabla_k\mathscr C_{ijl} =
-(b_{im} R_{jkl}{}^m+b_{jm} R_{kil}{}^m+b_{km} R_{ijl}{}^m). \label{CCCcomp}
\end{align}
Once again we read that a Codazzi tensor is Riemann compatible. By eq.\eqref{CCCcomp} 
the differential condition $\nabla_{(i}\mathscr C_{jk)l}=0$ is equivalent to the algebraic eq.\eqref{bR}.\\
A Veblen-like identity holds:
\begin{align}\label{VeblenC}
&\nabla_i \mathscr C_{jlk} + \nabla_j \mathscr C_{kil} + 
\nabla_k \mathscr C_{lji} + \nabla_l \mathscr C_{ikj}\\
&= b_{im} R_{jlk}{}^ m 
+ b_{jm} R_{kil}{}^m + b_{km} R_{lji}{}^m + b_{lm} R_{ikj}{}^m.\nonumber
\end{align}

\begin{example} For a concircular vector, $\nabla_i X_j =\rho g_{ij}$, the tensor $X_iX_j$ is Riemann compatible.\\
Proof: It is $\mathscr C_{jkl} =(\nabla_j \rho) g_{kl} - (\nabla_k \rho ) g_{jl}$ and $\nabla_i\mathscr C_{jkl} =(\nabla_i\nabla_j\rho) g_{kl} - (\nabla_i\nabla_k\rho ) g_{jl}$. The cyclic sum in \eqref{CCCcomp} gives zero. \\
Note: the existence of a concircular time-like vector is necessary and sufficient for a space-time to be 
generalized Robertson-Walker \cite{Chen_GRW}.
%In a Lorentzian manifold, the normalized vector $u_iu^i=-1$ is torse-forming: $\nabla_i u_j =\varphi (u_iu_j + g_{ij})$ (i.e. it is a %geodesic, shear-free and vorticity-free velocity field). If $\nabla_i\varphi =-u_i u^k\nabla_k\varphi $ then $u_iu_j$ is Riemann %compatible.
\end{example}
\begin{example}[{\bf Lovelock's identities}]{\quad}\\
{\bf 1)} The Codazzi deviation of the Ricci tensor is:
$ \mathscr C_{jkl} = \nabla_jR_{kl}-\nabla_k R_{jl} = -\nabla^m R_{jklm} $.
Property \eqref{CCCcomp} becomes a Lovelock's identity for the Riemann tensor {\rm (\cite{LovRun}, p.289)}: 
\begin{align}\label{lovelock}
\nabla_i\nabla^m R_{jklm} +\nabla_j\nabla^m R_{kilm} +\nabla_k\nabla^m 
R_{ijlm} = - R^m{}_{(i} R_{jk)lm}. 
\end{align}
{\bf 2)} The Codazzi deviation of Schouten's tensor\footnote{Schouten tensor: $S_{ij} =\frac{1}{n-2} \left [ R_{ij}-\frac{R}{2(n-1)} g_{ij}\right ] $. 
%$C_{jklm}= R_{jklm} +g_{jm}S_{kl}+g_{kl}S_{jm}-g_{jl}S_{km}-g_{km}S_{jl}$.
Properties: $\nabla_k S^k{}_j = \nabla_j S^k{}_k$,  $\nabla^m C_{jklm} = (n-3)(\nabla_k S_{jl}-\nabla_j S_{kl})$.}
 is $\mathscr C_{jkl} = -\frac{1}{n-3}\nabla^m C_{jklm}$. 
Property \eqref{CCCcomp} is $\nabla_{(i} \mathscr C_{jk) l}= -(n-3) S^m{}_{(i} R_{jk)lm}$. 
The term with the metric tensor in $S_{ij}$ does not contribute (Bianchi identity), and one is left with {\rm (see \cite{RCT})}:
\begin{align}
 \nabla_i\nabla^m C_{jklm} + \nabla_j\nabla^m C_{kilm} + \nabla_k\nabla^m C_{ijlm} =
- \frac{n-3}{n-2} R^m{}_{(i}R_{jk)lm} . \label{NNC}
\end{align}
In particular in $n>3$, if $\nabla_m C_{jkl}{}^m =0$ (conformally symmetric spaces, Roter \cite{Roter1972}) the Ricci tensor is Riemann compatible.
\end{example}
\begin{prop}
If $u_i u_j$ is Riemann compatible, and $u^ku_k\neq 0$, then $u_i$ is eigenvector of the Ricci tensor. 
\begin{proof}
Since $u_iu_j$ is Riemann compatible, it commutes with the Ricci tensor: $R_{ij}u^ju_k$ $= R_{kj}u^j u_i$.
Contraction with $u^k$ gives: $R_{ij}u^j (u_k u^k) = (R_{kj}u^j u^k ) u_i =0$. 
%In the other relation, $R_{ij}b^j{}_k - R_{kj}b^j{}_i =0$, contraction with 
\end{proof}
\end{prop}
We extrapolate a simple statement from Proposition 5.1 in \cite{Deszcz2013}. A direct proof is possible, by 
writing \eqref{bR} for the Ricci tensor in the warping coordinates:
\begin{prop}
In a warped spacetime $ds^2= \pm dt^2 + a(t)^2 g^*_{\mu\nu}dx^\mu dx^\nu $ the Ricci tensor is
Riemann compatible if and only if the Ricci tensor of the Riemannian submanifold  $(M^*,g^*)$ is 
compatible with the Riemann tensor of the submanifold:
$$ R^*_{\mu\sigma} R^*_{\nu\rho\lambda}{}^\sigma +R^*_{\nu\sigma} R^*_{\rho\mu\lambda}{}^\sigma +R^*_{\rho\sigma} R^*_{\mu\nu\lambda}{}^\sigma =0. $$
\end{prop}

\subsection{Geodesic maps} A map $(M,g)\to (M,\overline g)$  is 
{\em geodesic} if every geodesic line is mapped to a geodesic line. It is necessary and
sufficient that there exists a 1-form such that the Christoffel symbols are related by 
$\overline\Gamma_{ij}^k = \Gamma_{ij}^k +\delta_i{}^k X_j + X_i \delta^k{}_j$ 
(Levi-Civita, 1896).
%It implies $X_j = \frac{1}{n+1}(\overline\Gamma_{kj}^k - \Gamma_{kj}^k) $. By the relation 
%$\Gamma_{ki}^k=\frac{1}{2} \partial_i \log |g|$, it is
%$X_i = \nabla_i X$ with $X=\frac{1}{2(n+1)} \log |\overline g/g|$.
%By the relation $\partial_c g_{ab} = \Gamma^d_{ac} g_{bd} + \Gamma^d_{bc} g_{ad}$, the 
%Levi-Civita condition is equivalent to 
%$\nabla_k \overline g_{jl} = 2X_k\overline g_{jl} + X_j\overline g_{kl}+
%X_l\overline g_{kj}$.
The relation between the Riemann tensors is
\begin{align*}
\overline R_{jkl}{}^m = -\partial_j \overline \Gamma^m_{kl} +\partial_k \overline \Gamma^m_{jl} -
\overline\Gamma^d_{kl}\overline \Gamma^m_{jd} + \overline\Gamma^d_{jl} \overline\Gamma^m_{kd} 
= R_{jkl}{}^m  - \delta_k{}^m P_{jl} + \delta_j{}^m P_{kl},
\end{align*} 
where $P_{kl}=\nabla_kX_l-X_kX_l =P_{lk}$. 
It is: $\overline R_{jl} = R_{jl} + (n-1) P_{jl}$.\\
Geodesic maps preserve the (3,1) projective curvature tensor \cite{Sinyukov}: $\overline P_{jkl}{}^m=P_{jkl}{}^m$, where
$ P_{jkl}{}^m = R_{jkl}{}^m + \frac{1}{n-1} (\delta_j{}^m R_{kl} - \delta_k{}^m R_{jl} )$.
\begin{prop} [\cite{RCT}] 
If $b_{ij}=b_{ji}$, a geodesic map satisfies
\begin{equation}
b_{im} \overline R_{jkl}{}^m+b_{jm} \overline R_{kil}{}^m+
b_{km} \overline R_{ijl}{}^m =  
b_{im} R_{jkl}{}^m+b_{jm} R_{kil}{}^m+b_{km} R_{ijl}{}^m 
\end{equation}
Then, if $(R,b)$ is a compatible pair, also $(\overline R,b)$ is.
\end{prop}

\section{Weyl  compatible tensors}
A symmetric tensor  is Weyl compatible if:
\begin{align}
b_{im} C_{jkl}{}^m + b_{jm} C_{kil}{}^m + b_{km} C_{ijl}{}^m =0.
\end{align}
This identity holds for any symmetric tensor \cite{RCT}:
\begin{align}
b_{im} C_{jkl}{}^m + b_{jm} C_{kil}{}^m + b_{km} C_{ijl}{}^m 
= b_{im} R_{jkl}{}^m + b_{jm} R_{kil}{}^m + b_{km} R_{ijl}{}^m \label{RC}\\
+ \tfrac{1}{n-2}\left [ 
g_{kl} (b_{im} R_j{}^m - b_{jm} R_i{}^m ) + g_{il} (b_{jm} R_k{}^m - b_{km} 
R_j{}^m ) + g_{jl} (b_{km} R_i{}^m - b_{im} R_k{}^m )\right ].\nonumber
\end{align}
A simple consequence is obtained in dimension $n=3$, where the Weyl tensor is zero (see \cite{Deszcz2011}, in less simple manner):
 \begin{prop}
 In $n=3$ a Ricci tensor is Riemann compatible.
 \end{prop}
If $b_{ij}$ is Riemann compatible, then it commutes with the Ricci tensor.  As a result, the identity shows that $b_{ij}$ is also Weyl compatible. Therefore,
Riemann compatibility is a stronger condition than Weyl compatibility. The identity \eqref{RC} can be rewritten in terms of the Codazzi deviation:
\begin{align}
b_{im} C_{jkl}{}^m + b_{jm} C_{kil}{}^m + b_{km} C_{ijl}{}^m  =
 \nabla_i \mathscr D_{jkl} + \nabla_j \mathscr D_{kil} + \nabla_k \mathscr D_{ijl} \label {CCodD}\\
-\tfrac{1}{n-2} \nabla^m (\mathscr C_{ijm}g_{kl}+\mathscr C_{jkm}g_{il} +\mathscr C_{kim}g_{jl} ).  \nonumber   
%&\mathscr D_{jkl} = \mathscr C_{jkl} - \tfrac{1}{n-2}\left (\mathscr C_{jm}{}^m g_{kl}-\mathscr C_{km}{}^m g_{jl} \right ) \nonumber
\end{align}
where $\mathscr D_{jkl} = \mathscr C_{jkl} - \tfrac{1}{n-2}\left (\mathscr C_{jm}{}^m g_{kl}-\mathscr C_{km}{}^m g_{jl} \right )$.
\begin{example} If a vector field is torqued \cite{Chen}, i.e. $\nabla_i \tau_j = \rho g_{ij} + \alpha_i \tau_j$ with $\alpha_k \tau^k=0$, then 
$\tau_i \tau_j$ is Weyl compatible.\\
Proof:  one evaluates $\mathscr C_{jkl}=-\rho (\tau_j g_{kl}-\tau_k g_{jl})$ and $\mathscr D_{jkl} =-\frac{1}{n-2} \mathscr C_{jkl}$. It turns out that the r.h.s. of \eqref{CCodD} is zero.\\
Note: the existence of a torqued time-like vector is necessary and sufficient for a space-time to be twisted \cite{Chen}.
\end{example}
\begin{prop}[see remark 4.2 of \cite{Hervik}]
In a space-time of dimension $n=4$, if $u_iu_j$ is Weyl compatible and time-like unit ($u^ku_k=-1$) then the 
Weyl tensor is wholly determined by the electric tensor $E_{kl}=C_{jklm} u^ju^m $:
\begin{align}
C_{abcd} =  2(u_au_d E_{bc} -u_au_c E_{bd}+u_bu_c E_{ad}- u_bu_d E_{ac}) \label{CE}\\
+ g_{ad} E_{bc}- g_{ac} E_{bd}  + g_{bc} E_{ad} - g_{bd}E_{ac} \nonumber
\end{align}
\begin{proof}
In $n=4$ the following Lovelock's identity holds (\cite{LovRun}, ex 4.9 page 128):
\begin{align}
0= & g_{ar} C_{bcst} + g_{br} C_{cast}+ g_{cr} C_{abst} 
+g_{at} C_{bcrs} + g_{bt} C_{cars}+ g_{ct} C_{abrs} \nonumber \\
&+g_{as} C_{bctr} + g_{bs} C_{catr}+ g_{cs} C_{abtr} \nonumber
\end{align}
The contraction with $u^a u^r$ gives 
\begin{align}
0= & - C_{bcst} + u_b u^r C_{crst}+ u_c u^r C_{rbst} 
+u_t u^r C_{bcrs} + g_{bt} u^a u^r C_{cars}+ g_{ct} u^a u^r C_{abrs} \nonumber \\
&+u_s u^r C_{bctr} + g_{bs} u^a u^r C_{catr}+ g_{cs} u^a u^r C_{abtr} \nonumber\\
= & - C_{bcst} +  u^r (u_b C_{stcr}+ u_c  C_{rbst} 
+u_t C_{cbsr} + u_s C_{bctr} )\nonumber\\
&+g_{bt} E_{cs} - g_{ct} E_{bs} 
 - g_{bs} E_{ct}+ g_{cs} E_{bt} \nonumber
 \end{align}
 This gives the Weyl tensor in terms of its single and double contractions with $u^i$. 
If $u_iu_j$ is Weyl compatible, the single contraction is: $C_{jklr}u^r = u_k E_{jl}-u_j E_{kl}$, and the result is obtained. 
For an extension to $n>4$ see \cite{Hervik}.
\end{proof}
\end{prop}
%Space-times with Weyl-compatible tensor $k_ik_j$ where $k^ik_i=0$ (null vector) are Petrov type II or D.

\subsection{Conformal maps} A map $(M,g) \to (M,\hat g)$  is 
{\em conformal} if $\hat g_{kl} = e^{2\sigma} g_{kl}$.  The Christoffel symbols transform according to:
$\hat\Gamma^m_{ij} =\Gamma^m_{ij} + \delta^m{}_i X_j +X_i \delta^m{}_j - g_{ij} X^m, $
where $X_i =\nabla_i \sigma$. A conformal map leaves the Weyl tensor (3,1) unchanged: $\hat C_{jkl}{}^m =C_{jkl}{}^m$. 
Therefore, Weyl compatibility is an invariant property of conformal maps.

\section{K-compatible tensors}
Riemann and Weyl compatibility extend to $K-$compatibility, where $K$ is a generalised curvature tensor (GCT), i.e.
a tensor with the algebraic properties of the Riemann tensor under permutation of indices \cite{KobNom}:
%
%\begin{definition}  A tensor $K_{jklm}$ is a generalized curvature tensor {\rm (GCT)} if 
\begin{align}
&K_{jklm}=-K_{kjlm} = -K_{jkml} ,\label{K1}\\
&K_{jklm}+K_{kljm}+K_{ljkm}=0 ,\label{K2}\\
&K_{jklm}=K_{lmjk} \label{K3}.
\end{align}
%\end{definition}
%
\noindent
In analogy with the Riemann tensor, one shows that \eqref{K1} and \eqref{K2} imply the symmetry \eqref{K3}, and the
identity $K_{j(klm)}=0$. The tensor $K_{jl}=K_{jml}{}^m$ is symmetric. 

A symmetric tensor $b_{ij}$ is $K-$compatible if:
\begin{align}
b_i{}^m K_{jklm} + b_j{}^m K_{kilm} + b_k{}^m K_{ijlm}=0 \label{bKcomp}
\end{align}
and $(K,b)$ is a compatible pair. The property can be written $b^m{}_{(i} K_{jk)lm}=0$. \\
The metric tensor is $K-$compatible, by the Bianchi property \eqref{K2}. The tensors $b_{ij}$ and $K_{ij}$ commute:
$b_i{}^m K_{mk} -  K_{im} b^m{}_k=0$ (contract \eqref{bKcomp} with $g^{jl}$ and use symmetry).\\
Examples of $K-$compatible tensors were obtained by Shaikh et al. starting from specific metrics (see for example \cite{Shaikh1,Shaikh2}).
Bourguignon proved that if $b_{ij}$ is a Codazzi tensor then {\r R}$_{jklm} =R_{jkrs}b^r{}_l b^s{}_m$ is a GCT, \cite{Bourguignon}. 
We prove a more general statement:
\begin{prop}\label{pallino} If $a_{ij}$ and $b_{ij}$ are $K-$compatible, then 
$ \text{\r K}_{jklm}= K_{jkrs}(a^r{}_l b^s{}_m + b^r{}_l a^s{}_m) $ is a {\rm GCT}.
\begin{proof} The properties \eqref{K1} and \eqref{K3} are obvious; the Bianchi property \eqref{K2} completes the proof:
 \r K$_{(jkl)m}= a^r{}_{(l}{} K_{jk)rs} \,b^s{}_m + b^r{}_{(l}{} K_{jk)rs} \,a^s{}_m =0$ because each term is zero being 
 $a$ or $b$ $K-$compatible.
 \end{proof}
 \end{prop}
 
%\begin{example} 
%If the projector $h_{ij}$ on a submanifold $M^*$ is $K-$compatible, then $K^*_{mnpq}=K_{ijkl} h^i{}_m h^j{}_n h^k{}_p h^l{}_q$ is a %generalized curvature tensor in $M^*$.
%\end{example}

\subsection{Properties of $K-$compatible tensors}
A linear combination of K-compatible tensors obvioulsy is K-compatible. Now we prove: 
\begin{thrm} If $a$ and $b$ are $K-$compatible, then $\tfrac{1}{2}(ab+ba)$ is $K-$compatible.
\begin{proof}
Let $c{}_{ij} = a_i{}^k b_{kj}+
b_i{}^ka_{kj} $. Then:
\begin{align*}
&c^m{}_{(i} K_{jk)rm} = \, a_i{}^s b_s{}^m K_{jkrm} + a_j{}^s b_s{}^m K_{kirm} + a_k{}^s b_s{}^m K_{ijrm}
+ a\leftrightarrows b\\
& = -a_i{}^s (b_j{}^m K_{ksrm}+b_k{}^m K_{sjrm})-a_j{}^s (b_k{}^m K_{isrm}+b_i{}^m K_{skrm})\\
&\quad -a_k{}^s (b_i{}^m K_{jsrm}+b_j{}^m K_{sirm})+ a\leftrightarrows b\\
&= -(a_i{}^s b_j{}^m -a_j{}^s b_i{}^m) K_{ksrm} - (a_j{}^s b_k{}^m-a_k{}^s b_j{}^m)K_{isrm} \\
& \quad - (a_k{}^s b_i{}^m -a_i{}^s b_k{}^m )K_{jsrm}+   a\leftrightarrows b \\
&= -(a_i{}^s b_j{}^m -a_j{}^s b_i{}^m) (K_{ksrm} - K_{kmrs})- (a_j{}^s b_k{}^m-a_k{}^s b_j{}^m)(K_{isrm} - K_{imrs}) \\
& \quad - (a_k{}^s b_i{}^m -a_i{}^s b_k{}^m )(K_{jsrm}-K_{jmrs}) 
\end{align*}
\begin{align*}
&= (a_i{}^s b_j{}^m -a_j{}^s b_i{}^m) K_{krsm} +(a_j{}^s b_k{}^m-a_k{}^s b_j{}^m) K_{irsm} +
(a_k{}^s b_i{}^m -a_i{}^s b_k{}^m ) K_{jrsm}\\
&= (a_i{}^s b_j{}^m +b_i{}^s a_j{}^m) K_{krsm} +(a_j{}^s b_k{}^m+b_j{}^s a_k{}^m) K_{irsm} +
(a_k{}^s b_i{}^m +b_k{}^s a_i{}^m ) K_{jrsm}\\
&= \text{\r K}_{krij} + \text {\r K}_{irjk} + \text {\r K}_{jrki} = \text{\r K}_{(kri)j}=0
\end{align*}
because \r K is a GCT by Prop.\ref{pallino}. 
\end{proof}
\end{thrm}
Therefore, the linear space of $K-$compatible tensors is a special Jordan algebra.\\
In particular, the powers of $b$ are $K-$compatible (powers $n, n+1, ...$ are linear combinations of 
lower powers by Cayley-Hamilton theorem). 
In particular (with an exchange of indices) the tensor $(b^2)_j{}^s(b^2)_k{}^r K_{rslm} $ is a GCT. This enables the simple proof
 of the theorem in \cite{EXTDS}, so short that we reproduce it:

\begin{thrm}[Extended Derdzi\'nski-Shen theorem] Let $b_{ij}$ be K-compatible, $X^i$, $Y^i$, $Z^i$ be eigenvectors of 
$b_i{}^m$ with eigenvalues $x $, $y $, $z $.  If $x\neq z$ and $y\neq z$ then: 
\begin{align}
       K_{ijkl} X^i  Y^j Z^k =0. 
\end{align}
%(if $x\neq y\neq z$ any contraction of $K$ with the three vectors is zero; if two are equal the vanishing contraction has the two vectors %contracting with indices not linked by symmetry;  if the eigenvalues are equal nothing can be stated about vanishing contractions).
\begin{proof} 
Consider the identities $g^m{}_{(i}K_{jk)lm}=0$, $b^m{}_{(i}K_{jk)lm}=0$, $(b^2)^m{}_{(i}K_{jk)lm}=0$ and contract
them with $X^iY^jZ^k$. The three algebraic relations are
put in matrix form:
\begin{equation*}
\left[ \begin{array}{ccc}
1 & 1 & 1 \\
x & y & z \\
x^2 & y^2 & z^2 
\end{array}\right]
\left [ \begin{array}{c}
K_{jkli}X^iY^jZ^k\\
K_{kilj}X^iY^jZ^k\\
K_{ijlk}X^iY^jZ^k
\end{array}\right]=
\left [ \begin{array}{c} 0\\ 0\\ 0 \end{array}\right]
\end{equation*}
The determinant of the matrix is
$(x-y)(x-z)(z-y)$. If the eigenvalues are all different 
then $K_{ijkl}X^iY^jZ^k=0$ (with contraction of any three indices). 
If $x=y\neq z$, the reduced system of equations still implies $K_{ijkl}X^iY^jZ^k=0$. 
\end{proof}
\end{thrm}
\begin{prop} If $b$ is K-compatible and invertible, then $b^{-1}$ is K-compatible:
\begin{align}
(b^{-1})^j{}_{(s} K_{r l) kj} =0
\end{align}
\begin{proof} Multiply \eqref{bKcomp} by $(b^{-1})^i{}_r (b^{-1})^j{}_s $ and obtain the identity:
$ (b^{-1})^j{}_s K_{jklr} + (b^{-1})^i{}_r K_{kils} + (b^{-1})^i{}_r (b^{-1})^j{}_s b^m{}_k K_{ijlm}=0$. 
Rewrite it as: 
$$ (b^{-1})^j{}_{(s} K_{r l) kj} - (b^{-1})^j{}_l K_{sr kj} + (b^{-1})^i{}_r (b^{-1})^j{}_s b^m{}_k K_{ijlm}=0$$
The last two terms cancel, as shown by the chain:\\
$(b^{-1})^j{}_l K_{sr kj} = (b^{-1})^i{}_r (b^{-1})^j{}_s b^m{}_k K_{ijlm}$ $\Leftrightarrow$
$K_{sr kb} b^r{}_a = b^i{}_b (b^{-1})^j{}_s b^m{}_k K_{a j lm}$ \\ $\Leftrightarrow$
$b^s{}_c K_{sr kb} b^r{}_a = b^l{}_b  b^m{}_k K_{ac lm}$ $\Leftrightarrow$
\r K$_{kbc a}=$\r K$_{ac b k}$, which is true as \r K is a GCT.
\end{proof}
\end{prop}

We prove a Veblen-like identity:
\begin{prop} If $b_{ij}$ is K-compatible then:
\begin{align}
b_i{}^m K_{jklm} - b_j{}^m K_{ilkm} + b_k{}^m K_{iljm} - b_l{}^m K_{jkim}=0.
\end{align}
\begin{proof} $ 0 = b_i{}^m K_{jklm} + b_j{}^m K_{kilm} + b_k{}^m K_{ijlm}=
b_i{}^m K_{jklm} - b_j{}^m (K_{ilkm}+K_{lkim}) + b_k{}^m K_{ijlm} 
=b_i{}^m K_{jklm} - b_j{}^m K_{ilkm} +b_l{}^m K_{kjim}+b_k{}^m K_{jlim} + b_k{}^m K_{ijlm} \\
=b_i{}^m K_{jklm} - b_j{}^m K_{ilkm} +b_l{}^m K_{kjim}-b_k{}^m K_{lijm} $.
\end{proof}
\end{prop}

\subsection{More on generalised curvature tensors} A linear combination of GCTs is a GCT. 
Given two compatible pairs $(K,a)$ and $(K,b)$ a new GCT tensor is obtained in Prop.\ref{pallino}. 
In particular, if $a_{ij}=g_{ij}$ (the metric tensor) the following $K'$ is a {\rm GCT}:
\begin{align}
K'_{jklm}=K_{jkrs} (\delta^r{}_l b^s{}_m + b^r{}_l \delta^s{}_m)=K_{jkls} b^s{}_m-K_{jkms} b^s{}_l
\end{align}

\begin{prop} 
If $b$ is K-compatible, then $b$ is $K^\prime$-compatible.
\begin{proof} The tensor $K^\prime_{jklm}=K_{jklr}b^r{}_m - K_{jkmr} b^r{}_l$ is a GCT.  Let us evaluate:
$b^m{}_i K^\prime_{jklm} = b^m{}_{i}K_{jklr} b^r{}_m- b^m{}_{i} K_{jkmr}b^r{}_l = (b^2)^r{}_{i}K_{jklr} - \text{\r K}_{jkim}$. Both tensors vanish if the cyclic sum $(ijk)$ is taken. 
\end{proof}
\end{prop}
\begin{prop}
$(K,b)$ is a compatible pair for any symmetric tensor $b$ if and only if
\begin{align}
K_{ijlm} = \frac{K}{n(n-1)}(g_{il}g_{jm}-g_{im}g_{jl}) \label{Kgg}
 \end{align}
where $K$ is a scalar field.
\begin{proof}
The symmetry of the tensor is made explicit by writing $b_{ij} =\frac{1}{2}  b^{rs} (g_{ir}g_{js} + g_{is}g_{jr}) $. The compatibility relation 
must hold for any $b^{rs}$, then: 
$$0=g_{ir}K_{jkls} + g_{jr} K_{kils} + g_{kr} K_{ijls} + g_{is}K_{jklr} + g_{js} K_{kilr} + g_{ks} K_{ijlr}. $$
Contraction with $g^{ks}$ gives $ (n-1) K_{ijlr} = g_{jr}K_{il}-g_{ir}K_{jl}$; contraction with $g^{il}$ gives $K_{jr}=\frac{1}{n} g_{jr} K^i{}_i$ and \eqref{Kgg} follows.
The reverse, i.e. \eqref{Kgg} implies \eqref{bKcomp}, is shown by direct check. 
\end{proof}
\end{prop}
 
A pseudo-Riemannian manifold of dimension $n>2$ is an {\em Einstein manifold} if $R_{ij} =\frac{1}{n}R g_{ij}$ where $R$ is the scalar curvature. Since $\nabla_iR^i{}_j = \frac{1}{2}\nabla_j R$, the scalar curvature is constant. 
A manifold is a {\em constant curvature manifold} if the Riemann tensor has the form 
\eqref{Kgg}. Such manifolds are Einstein manifolds.
\begin{cor} A manifold is a constant curvature manifold if and only if\\ 
 $b_i{}^m R_{jklm} + b_j{}^m R_{kilm} + b_k{}^m R_{ijlm}=0$ for all symmetric tensors. 
 \end{cor}
 % 

% \begin{cor}
%A conformally flat Einstein manifold is a constant curvature manifold.
% \begin{proof}
% The relation \eqref{RC} simplifies. For Einstein manifolds symmetric tensors commute
% with the Ricci tensors, then $b_i{}^m R_{jklm} + b_j{}^m R_{kilm} + b_k{}^m R_{ijlm}=0$ 
% for all symmetric tensors, i.e. the manifold is constant curvature.
 % \end{proof}
 %\end{cor}

%Given the compatible pair $(K,b)$, the tensor $b_{ij}$ commutes with 
%$K_{ij}$ and with \r b$_{ij}$. By the Jordan identity, $b_{ij}$ commutes with
%the antisymmetric tensor $[K , \text {\r b}]$.

\end{document}